\documentclass[12 pt]{amsart}

\usepackage{graphicx}
\usepackage{color}
\usepackage{amssymb, amsmath, amsthm}
\usepackage{mathptmx}

\usepackage{epsfig}
\usepackage{epstopdf}
\usepackage{graphicx}

\newtheorem{thm}{Theorem}[section]

\newtheorem{cor}[thm]{Corollary}

\newtheorem{defn}[thm]{Definition}
\newtheorem{lemma}[thm]{Lemma}

\newcommand{\R}{\mathbb R} %REALS

\newcommand{\bi}{\begin{itemize}}
\newcommand{\ei}{\end{itemize}}
\newcommand{\be}{\begin{enumerate}}
\newcommand{\ee}{\end{enumerate}}

\newcommand{\n}{\beta}

\newcommand{\emp}{\emptyset}
\newcommand{\X}{\times}

\newcommand{\A}{\alpha}

\newcommand{\K}{\mathcal{K}}

\begin{document}
\title{Unexpected local minima in the width complexes for knots}
\author{Alexander Zupan}

\maketitle

\begin{abstract}
In \cite{complex}, Schultens defines the width complex for a knot in order to understand the different positions a knot can occupy in $S^3$ and the isotopies between these positions.  She poses several questions about these width complexes; in particular, she asks whether the width complex for a knot can have local minima that are not global minima.  In this paper, we find an embedding of the unknot $0_1$ that is a local minimum but not a global minimum in the width complex for $0_1$.  We use this embedding to exhibit for any knot $K$ infinitely many distinct local minima that are not global minima of the width complex for $K$.
\end{abstract}

\section{Introduction}
In \cite{gabai}, Gabai defines knot width and thin position as a measure of the complexity of various Morse functions on a given knot in $S^3$.  One important aspect of thin position is that it yields an embedding of a given knot that is minimal with respect to certain types of isotopies.  In \cite{complex}, Schultens defines the width complex of a knot in order to better understand these isotopies and the various positions a given knot can occupy in $S^3$.  Specifically, she asks the following two questions: \\

\noindent \textbf{Question 12.} \emph{
Can the width complex of a knot have local minima that are not global minima?} \\

\noindent \textbf{Question 13.} \emph{Is every vertex of the width complex of a knot connected to one of the global minima of this complex by a monotonically decreasing path?} \\

Schultens also defines a similar width complex for 3-manifolds, and her Theorem 13 from \cite{complex} provides a positive answer to the 3-manifold version of Question 12, namely that there exist 3-manifolds whose width complexes contain local minima that are not global minima.  On the other hand, combining the results of \cite{bonotal}, \cite{schar}, and \cite{wald}, we see that if $M$ is $S^3$ or a lens space, then the width complex of $M$ has a unique minimum, corresponding to a minimal genus Heegaard splitting.  Thus, it seems reasonable to expect that the simplest knots might share this property.  This is further suggested by Otal's proof that non-minimal bridge positions of the unknot and 2-bridge knots are stabilized \cite{otal} and Ozawa's recent proof of the same statement for torus knots \cite{ozawa}. \\

Schultens compares Question 13 to one answered by Goeritz in 1934.  Goeritz produced a nontrivial diagram of the unknot $0_1$ such that any Reidemeister move increases the diagram's crossing number.  We find an analogous result concerning the width complex of $0_1$, finding a nontrivial embedding such that any isotopy must increase the complexity of the embedding, the difference being that in this context complexity refers to knot width instead of crossing number.  As a result, we give an affirmative answer to Schultens' first question, which shows that the answer to the second question must be no.  In fact, we show the surprising and much stronger result that for every knot $K$, the width complex of $K$ has infinitely many local minima that are not global minima. \\

\section{Definitions}
Let $K$ be a knot in $S^3$, and fix a Morse function $h:S^3 \rightarrow \R$ such that $h$ has exactly two critical points.  We can think of $K$ as an equivalence class, denoted $\K$, of the set of embeddings of $S^1$ into $S^3$ modulo ambient isotopy.  In the usual definition of knot width, the embedding of $K$ is fixed and the Morse function $h$ is allowed to vary up to isotopy; however, this definition is equivalent with the one that follows.  Let $k \in \K$ such that $h \mid_k$ is Morse, and let $c_0 < c_1 < \dots < c_n$ be the critical levels of $h\mid_k$.  Choose regular levels $c_0 < r_1 < c_1 < \dots < r_n < c_n$, and define
\[ w(k) = \sum_{i=1}^n |h^{-1}(r_i) \cap K| .\]
Now, let
\[ w(K) = \min_{k \in \K} w(k).\]
The invariant $w(K)$ is called the width of $K$, and if $k \in \K$ satisfies $w(K) = w(k)$, we say that $k$ is a thin position for $K$. \\

For our purposes it will be useful to split an embedding $k$ into thick and thin levels.  For $2 \leq i \leq n-1$, we say that a regular value $r_i$ of $h \mid_k$ corresponds to a thick level $R_i = h^{-1}(r_i)$ if $|h^{-1}(r_i) \cap k| > |h^{-1}(r_{i-1}) \cap k|,|h^{-1}(r_{i+1}) \cap k|$.  Likewise, $r_i$ corresponds to a thin level $R_i = h^{-1}(r_i)$ if $|h^{-1}(r_i) \cap k| < |h^{-1}(r_{i-1}) \cap k|,|h^{-1}(r_{i+1}) \cap k|$.  Let $a_0,\dots,a_m$ ($b_1,\dots,b_m$) represent the regular values of $h \mid_k$ corresponding to thick (thin) levels $A_0,$ $\dots,A_m$ ($B_1,\dots,B_m$), where $a_0 < b_1 < a_1 < \dots < b_m < a_m$. \\

Note that $h^{-1}([b_i,a_i]) \cap k$ consists of vertical segments and arcs $\A_1,\dots,\A_l$, $l\geq 1$, where each $\A_j$ has exactly one minimum and is isotopic to an arc $\n_j$ in $A_i$.  In this case, $\A_j$ cobounds a disk $D$ with $\n_j$ such that $D$ has no critical points with respect to $h$ in its interior.  We call $D$ a \emph{strict lower disk for $k$ at $A_i$}. For any $r < c_0$, the lowest minimum of $h\mid_k$, we have that $h^{-1}([r,a_1]) \cap k$ consists of arcs $\A_1,\dots,\A_l$, which cobound pairwise disjoint strict lower disks with arcs $\n_1,\dots,\n_l$ contained in $A_1$.  \\

Similarly,  $h^{-1}([a_i,b_{i+1}]) \cap k$ consists of vertical segments and arcs $\A_1,\dots,$ $\A_l$, $l\geq 1$, where each $\A_j$ has exactly one maximum and is isotopic to an arc $\n_j$ in $A_i$.  Here $\A_j$ cobounds a disk $E$ with $\n_j$ such that $E$ has no critical points in its interior, and we call $E$ a \emph{strict upper disk for $k$ at $A_i$}. For any $r > c_n$, the highest maximum of $h\mid_k$, we have that $h^{-1}([a_n,r]) \cap k$ consists of arcs $\A_1,\dots,\A_l$, which cobound pairwise disjoint strict upper disks with arcs $\n_1,\dots,\n_l$ contained in $A_n$.  \\

Let $k,k' \in K$ with corresponding thick/thin levels $A_0,B_1,A_1,\dots,B_l,A_l$ and $A'_0,B'_1,\dots,B'_{l'},A'_{l'}$.  We say that $k \sim k'$ if $l = l'$ and there is an isotopy of $S^3$ taking $k$ to $k'$, $A_i$ to $A'_i$, and $B_i$ to $B'_i$.  In this case, we call this isotopy a \emph{level isotopy}, and we have $w(k) = w(k')$, so that $k$ and $k'$ carry exactly the same information with respect to width and to upper and lower disks.  Thus, from this point forward we will (under slight abuse of notation) let $\mathcal{K}$ denote the set of embeddings isotopic to $K$ up to this equivalence. \\

\section{The Width Complex of $K$}

Now, we use the set $\K$ and pairs of strict upper and lower disks to define the width complex of $K$, a directed graph $\Gamma$ whose vertices correspond to elements of $\mathcal{K}$.  We first make several definitions:
\begin{defn}
Suppose that $k \in \K$.  If $(D,E)$ is a pair of strict upper and lower disks for a thick level $A_i$ such that $D \cap E$ is a single point in $k$, we say that $A_i$ is \emph{stabilized}.  If $D \cap E = \emp$, we say that $A_i$ is \emph{weakly reducible}.  In either case, we say that $A_i$ is \emph{reducible}.  If $A_i$ is not reducible, then $A_i$ is \emph{strongly irreducible}.
\end{defn}
Elements of $k \in \K$ with reducible thick surfaces will be at the tail of directed edges in the width complex of $K$.  If $k \in \K$ has a stabilized thick surface $A_i$, we can slide $k$ along the pair $(D,E)$ of upper and lower disks for $A_i$ to cancel out a minimum and maximum, changing $k$ to $k' \in \K$ such that $w(k') = w(k) - (2|A_i \cap k| - 2)$.  As in \cite{complex}, we call this a \emph{Type I move}.  If $k$ has a weakly reducible thick surface $A_i$, we can again slide $k$ along the pair $(D,E)$ to move a minimum of $k$ above a maximum of $k$.  This changes $k$ to $k' \in \K$ such that $w(k') = w(k) - 4$, and we call this a \emph{Type II move}.  In either case, we call $(D,E)$ a \emph{pair of reducing disks} at $A_i$ and we make an directed edge from  $k$ to $k'$ in $\Gamma$.  The next theorem, Theorem 1 from \cite{complex}, is important to our understanding of the width complex:
\begin{thm}
The width complex of a knot is connected.
\end{thm}
This theorem says that given $k,k' \in \K$, there is a series of level isotopies and Type I and Type II moves taking $k$ to $k'$.  Schulten's width complex also contains higher dimensional cells, but we need only consider the one-skeleton of the complex in this context.
\begin{defn}
We call $k \in \K$ a \emph{local minimum} of the width complex if there are no directed edges leaving $k$ in $\Gamma$.
\end{defn}
The position $k$ is called a local minimum because any isotopy that changes $k$ to $k' \in \K$ must increase $w(K)$.  Let $\hat{\K} \subset \K$ denote the set of local minima of the width complex of $K$.  It is clear that any thin position $k$ for $K$ must come from $\hat{\K}$; otherwise there is an isotopy decreasing $w(k)$.  We also have the following, the proof of which is clear from the definition of the width complex:
\begin{lemma}\label{mini}
An element $k \in \K$ is in $\hat{\K}$ if and only if every thick level of $k$ is strongly irreducible.
\end{lemma}
Using the definitions of this section, we can reformulate Schulten's questions as follows: \\

\noindent \textbf{Question 12.} \emph{
Is there a knot $K$ with $k \in \hat{\K}$ such that $w(k) > w(K)$?} \\

\noindent \textbf{Question 13.} \emph{
Given $k \in \K$, is there a directed path in $\Gamma$ starting at $k$ and ending at a thin position for $K$?} \\

Explicitly, a directed path is a sequence of vertices $k=k_0,k_1,\dots,k_n$ such that each there is a directed edge from $k_i$ to $k_{i+1}$ for each $i <n$.

\section{A local minimum in the width complex of the unknot}

Let $K$ be the unknot in $S^3$, and let $k \in \K$ be the position of the unknot depicted in Figure 1, where $h$ is the standard height projection onto a vertical axis.
\begin{figure}\label{unkn}
  \centering
      \includegraphics[width=0.55\textwidth]{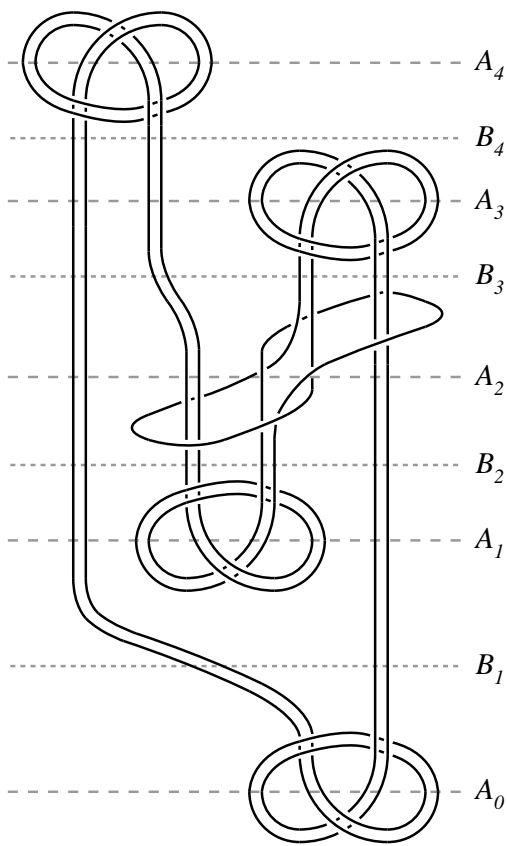}
  \caption{A troublesome embedding $k$ of the unknot, shown with thick/thin levels}
\end{figure}

We will label the thick/thin levels of $k$ as $A_0,B_1,A_1,B_2,A_2,B_3,A_3,B_4,A_4$, as shown.  First, we need several results about bridge position.
\begin{defn}
For any knot $K$ with embedding $k$, the bridge number of $k$, $b(k)$, is defined to be the number of maxima in $k$ with respect to $h$, and the bridge number of $K$, $b(K)$, is the minimum of $b(k)$ over $k \in \K$.  We call $k$ a \emph{bridge position} for $K$ if $b(k) = b(K)$ and $k$ has exactly one thick level, called a bridge sphere for $k$.
\end{defn}

Schultens shows in \cite{bridge} that the bridge number of any $(n,\ast)$-cable of a 2-bridge knot is $2n$, and we demonstrate in \cite{zupan} that any thin position is a bridge position for such a knot.  In this case, the bridge sphere must be strongly irreducible.  We will use this fact in the following:
\begin{thm}\label{unknot}
The pictured embedding $k$ of the unknot is a local minimum in the width complex.
\begin{proof}
By Lemma \ref{mini}, it suffices to show that every thick surface of $k$ is strongly irreducible.  Observe that the thick surfaces $A_0,A_1,A_3,A_4$ are identical except for the extra vertical segments contained in $A_1$ and $A_3$.  Thus, we need only show that $A_0$ and $A_2$ are strongly irreducible. \\

\textbf{Claim 1}: $A_0$ is strongly irreducible.  Suppose not.  Then there is a pair of reducing disks $(D,E)$ at $A_0$.  Let $b_1$ denote the regular value corresponding to the thin level $B_1$.  If we restrict our attention to $k^*  = k \cap h^{-1}(-\infty,b_1]$, we can easily see that by adding two arcs to the four intersection points of $k$ with $B_1$, we can complete $k^*$ to a $(2,\ast)$-cable of the trefoil, whose thin position is bridge position by the discussion above, and such that $A_0$ becomes a bridge sphere.  Thus, the pair $(D,E)$ of reducing disks at $A_0$ is also a pair of reducing disks at the bridge sphere $A_0$ of the trefoil's cable, a contradiction to the fact that this cable is in thin position.  We conclude that $A_0$ and thus $A_1$, $A_3$, and $A_4$ are strongly irreducible. \\

\begin{figure}\label{thksurf}
  \centering
      \includegraphics[width=0.6\textwidth]{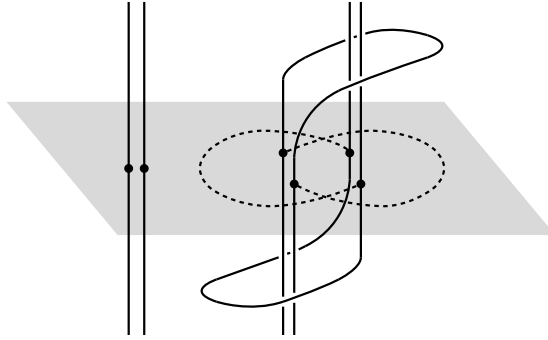}
  \caption{The set $k' \subset k$ contained in $\mathcal{A}_2$, with projections of upper and lower disks onto $A_2$.}
\end{figure}

\textbf{Claim 2}: $A_2$ is strongly irrreducible.  Let $b_2$ and $b_3$ be the critical values corresponding to $B_2$ and $B_3$, respectively.  Then $\mathcal{A}_2 = h^{-1}([b_2,b_3])$ is homeomorphic to $S^2 \X I$, and $k$ has exactly one minimum contained in an arc $\kappa_1$ and exactly one maximum contained in an arc $\kappa_2$ properly embedded in $\mathcal{A}_2$.  If we push $\kappa_1$ and $\kappa_2$ along an obvious pair of lower and upper disks $(D,E)$ onto arcs in $A_2$, we see the picture in Figure 2.  Note that $\mathcal{A}_2$ contains six additional vertical segments, four of which intersect $D$ or $E$ and two of which do not, call these $\gamma_1$ and $\gamma_2$.  Let $k' \subset \mathcal{A}_2 = \kappa_1 \cup \kappa_2 \cup \gamma_1 \cup \gamma_2$, shown in Figure 2. \\

Suppose by way of contradiction that $A_2$ is reducible.  Then there is a pair of reducing disks $(D',E')$ for $k$ contained in $\mathcal{A}_2$, and certainly $(D',E')$ also constitutes a pair of reducing disks for $k'$ in $\mathcal{A}_2$.  We note that $D' \cap E' = \emp$, since $\kappa_1$ and $\kappa_2$ do not share endpoints.  Consider the link $L$ shown on the left of Figure 3, with standard height projection $p$.  There are regular values $x_1$ and $x_2$ such that the pair $(p^{-1}([x_1,x_2]),L \cap p^{-1}([x_1,x_2]))$ is homeomorphic to $(\mathcal{A}_2,k')$.  But this implies that $L$ is isotopic to the unlink by an isotopy along $D'$ and $E'$, the result of which is shown on the right of Figure 3.  This contradicts the fact that $L$ is not the unlink (see for instance, \cite{rolf}, Section 5E).  We conclude that $A_2$ is strongly irreducible, completing the proof.

\begin{figure}\label{linkage}
  \centering
      \includegraphics[width=0.8\textwidth]{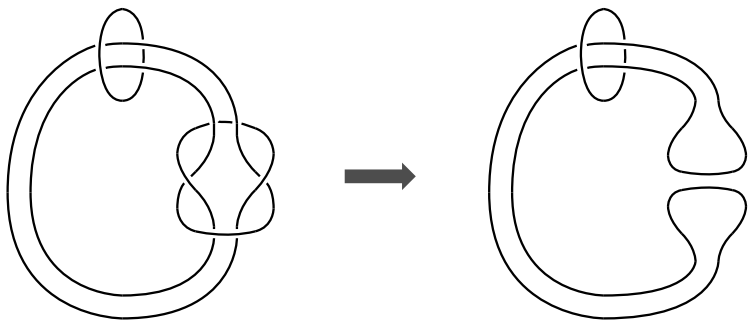}
  \caption{If $A_2$ is reducible, the above link must be isotopic to the unlink}
\end{figure}

\end{proof}
\end{thm}

\section{Local minima in the width complex of an arbitrary knot}
Suppose that $k_1$ and $k_2$ are embeddings representing local minima in the width complexes of knots $K_1$ and $K_2$.  Then we can find an embedding $k$ of $K_1 \# K_2$ by connecting the highest maximum of $k_1$ to the lowest minimum of $k_2$.  Observe that this creates a new thin surface but does not interfere with the reducibility of the thick surfaces of $k_1$ and $k_2$.  Thus, every thick surface of $k$ is strongly irreducible, and by Lemma \ref{mini}, $k$ represents a local minimum in the width complex of $K_1 \# K_2$.  For instance, consider the projection of the figure eight $4_1$ knot shown in Figure 4.  Note that bridge position is thin position for $4_1$ by \cite{thomps}.  Here we have taken $k_1$ to be bridge position of the figure eight knot and $k_2$ to be the unknot projection shown above, creating a new projection $k$ of $4_1$.  Since every thick sphere is strongly irreducible, this projection is a local minimum in the knot's width complex.  This suggests the following:

\begin{figure}\label{fig8}
  \centering
      \includegraphics[width=0.5\textwidth]{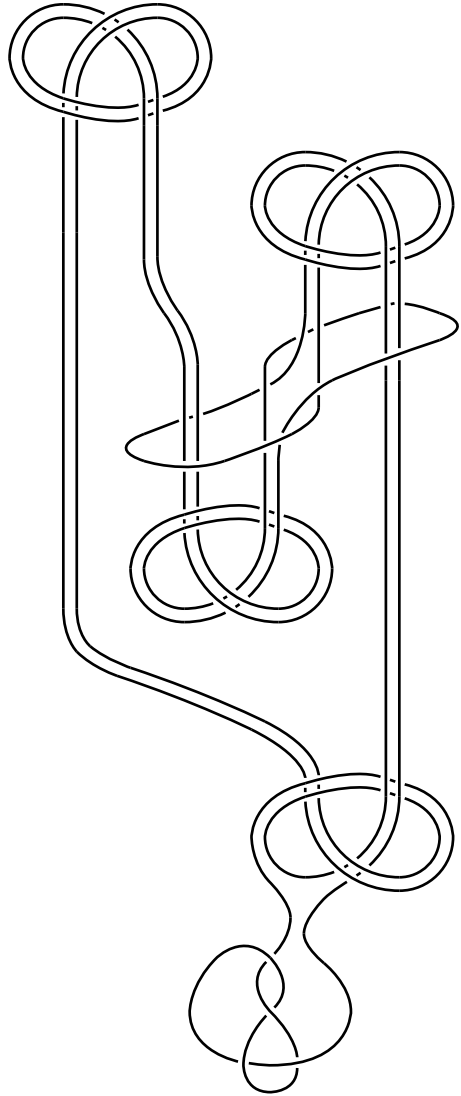}
  \caption{A local minimum in the width complex of the figure eight knot}
\end{figure}
\begin{cor}
The width complex of every knot contains infinitely many local minima.
\begin{proof}
Let $K$ be an arbitrary knot, with embedding $k$ representing a local minima in the width complex of $K$.  For any such $k$, we exhibit another local minima $k'$ of the width complex of $K$ with $w(k') > w(k)$, showing that there are infinitely many such embeddings.  Let $K_0$ denote the unknot, and let $k_0$ be the embedding representing the local minimum of the width complex in Theorem \ref{unknot}.  Since $K \# K_0 = K$, we can attach $k$ to $k_0$ by connecting the highest maximum of $k$ to the lowest minimum of $k_0$ to get a new embedding $k'$ of $K$ with $w(k') > w(k)$.  By the above argument, every thick sphere of $k'$ is strongly irreducible, so $k'$ is another local minimum in the width complex.
\end{proof}
\end{cor}


\begin{thebibliography}{5}
\bibitem{ozawa} M Ozawa, Non-minimal bridge positions of torus knots are stabilized, arxiv:1006.1026.
\end{thebibliography}
\end{document}